\documentclass[preprint,authoryear,11pt]{elsarticle}
\usepackage{amssymb,amsmath,latexsym}
\usepackage[T1]{fontenc}
\usepackage{amsmath,amsfonts,amssymb}
\usepackage[english]{babel}
\usepackage{latexsym}
\usepackage{epsfig}
\usepackage{graphics}
\usepackage{enumerate}
\usepackage{multirow}

\usepackage{tabularx}
\usepackage[left=2.5cm,right=2cm,top=2cm,bottom=2cm]{geometry}

 \usepackage{natbib}
\numberwithin{propo}{section}
\newtheorem{thm}{Theorem}
\numberwithin{thm}{section}

\numberwithin{rem}{section}

\numberwithin{exemple}{section}

\numberwithin{lem}{section}

\def\zak{\null\hfill{$\Box$}\par\vspace*{0.2cm}}

\begin{document}

\begin{frontmatter}

\title{Portmanteau test for the asymmetric power GARCH model when the power is unknown}

\author[Yacouba]{Y. Boubacar Ma{\"\i}nassara}
\address[Yacouba]{Universit\'e Bourgogne Franche-Comt\'e, \\
Laboratoire de math\'{e}matiques de Besan\c{c}on, \\ UMR CNRS 6623, \\
16 route de Gray, \\ 25030 Besan\c{c}on, France.}
\ead{mailto:yacouba.boubacar\_mainassara@univ-fcomte.fr}

\author[Yacouba]{O. Kadmiri}
\ead{mailto:othman.kadmiri@univ-fcomte.fr}

\author[Yacouba]{B. Saussereau}
\ead{mailto:bruno.saussereau@univ-fcomte.fr}

\begin{abstract}
It is now widely accepted that, to model the dynamics of daily financial returns, volatility models have to incorporate the so-called leverage effect.
 We derive the asymptotic behaviour of the squared residuals  autocovariances  for the class of asymmetric power GARCH model when the power is unknown and is jointly estimated with the model's parameters. We then deduce a portmanteau adequacy test based on the autocovariances of the squared residuals. These asymptotic results are illustrated by Monte Carlo experiments. An application to real financial data is also proposed.
\end{abstract}

\end{frontmatter}

\section{Introduction}\label{Introduction}
The autoregressive conditional heteroscedastic (ARCH) model introduced by \cite{E-ARCH} expresses the conditional variance (volatility) of the process as a
linear functional of the squared past values. This model has a lot of extensions. For instance, \cite{B-GARCH} generalized the ARCH (GARCH) model by adding the past realizations of the volatility.
The GARCH models are also characterized by a volatility specified as a linear function
of the squared past innovations.
Thus, by construction, the conditional variance only depends on the modulus of the
past variables: past positive and negative innovations have the same effect on the current volatility.
This property is in contradiction with many empirical studies on series of stocks, showing a negative
correlation between the squared current innovation and the past innovations. For instance,
%The weakness of this model is the interpretation of the past return because it does not make the difference between the past positive returns and the past negative returns. In 1976, Black\nocite{B}
\cite{B} showed that the past negative returns seem to have more impact on the current volatility than the past positive returns.
Numerous financial series present this stylised fact, known as the leverage effect.
%The so-called leverage effect was noted by Black (1976), and involves
%an asymmetry of the impact of past positive and negative values on the current volatility.
Since 1993, a lot of extensions are made to consider the leverage effect. Among the various asymmetric GARCH processes
introduced in the econometric literature, the more general is the asymmetric power GARCH (APGARCH for short) model of \cite{DGE}. For some positive constant $\delta$, it is defined by
\begin{equation}\label{model}
\left\{\begin{aligned}
\varepsilon_t &= \zeta_t\eta_t\\
\zeta_t^\delta &= \omega_0 + \sum\limits_{i=1}^q \alpha_{0i}^+ (\varepsilon_{t-i}^+)^\delta + \alpha_{0i}^- (-\varepsilon_{t-i}^-)^\delta + \sum\limits_{j=1}^p \beta_{0j} \zeta_{t-j}^\delta,
\end{aligned}\right.
\end{equation}
where $x^+ = \max(0,x)$ and $x^- = \min(0,x)$. It is assumed that

	\indent \textbf{A0:} $(\eta_t)$ is a sequence of independent and identically distributed (iid, for short) random variables with $\mathbb{E}\vert \eta_t\vert^r< \infty$ for some $r> 0$.

In the sequel, the vector of parameter of interest (the true parameter) is denoted $\vartheta_0= (\omega_0, \alpha_{01}^+, \ldots, \alpha_{0q}^+, \alpha_{01}^-, \ldots, \alpha_{0q}^-, \beta_{01},\ldots, \beta_{0p},\delta)'$ and satisfies the positivity constraints  $\vartheta_0\in ]0,+\infty[ \times [0,+\infty[^{2q+p}\times ]0,+\infty[$.
%We write the vector of parameters $\vartheta := (\theta', \tau)' = (\omega, \alpha_1^+, \ldots, \alpha_q^+, \alpha_1^-, \ldots, \alpha_q^-, \beta_1,\ldots, \beta_p,\tau)'$. The parameter of interest (the true parameter) is denoted $\vartheta_0 :=(\theta_0',\delta)'$, where $\theta_0= (\omega_0, \alpha_{01}^+, \ldots, \alpha_{0q}^+, \alpha_{01}^-, \ldots, \alpha_{0q}^-, \beta_{01},\ldots, \beta_{0p})'$.
The representation \eqref{model} includes various GARCH time series models: the standard
GARCH of \cite{E-ARCH} and \cite{B-GARCH} obtained for $\delta=2$ and $\alpha_{0i}^+=\alpha_{0i}^-$ for $i=1,\dots,q$;
the threshold ARCH (TARCH) model of \cite{RZ-TGARCH} for $\delta = 1$ and the GJR model of \cite{GJR-GJR} for $\delta= 2$.

After identification and estimation of the GARCH processes, the next important step in the GARCH modelling consists
in checking if the estimated model fits the data satisfactorily. This adequacy checking step allows to validate or invalidate the choice of the orders $p$ and $q$.
Thus it is important to check the validity of a GARCH$(p,q)$ model, for given orders $p$ and $q$. This paper is devoted to the
problem of the validation step of APGARCH$(p,q)$ representations \eqref{model} processes, when the power $\delta$ is estimated.
Based on the residual empirical autocorrelations,  \cite{BP} derived a goodness-of-fit test, the portmanteau test, for univariate
strong autoregressive moving-average (ARMA) models (i.e. under the assumption that the error term is iid).  \cite{LB} proposed a modified portmanteau test which is nowadays one
of the most popular diagnostic checking tool in ARMA modelling of time series.
Since the articles by \cite{LB} and \cite{McL}, portmanteau tests have been important tools in time series analysis, in particular
for testing the adequacy of an estimated ARMA$(p,q)$ model. See also \cite{Li-Book}, for a reference book on the portmanteau tests.

The intuition behind these portmanteau tests is that if a given
time series model with iid innovation $\eta_t$ is appropriate for the data at hand,
the autocorrelations of the residuals $\hat{\eta}_t$ should be close to zero, which is the theoretical
value of the autocorrelations of  $\eta_t$. The standard portmanteau tests thus consists in rejecting
the adequacy of the model for large values of some quadratic form of the residual
autocorrelations.

\cite{LM} and \cite{LL} studied a portmanteau test based on the autocorrelations of the squared residuals. Indeed the test based on the autocorrelations is irrelevant because the process such that this use to define a GARCH model ($\hat{\eta}_t = \varepsilon_t/ \hat{\sigma}_t$) with $\hat{\sigma}_t$ independent of $\sigma\{\eta_u, u < t\}$, is a martingale difference and thus is uncorrelated. Concerning the GARCH class model, \cite{BHK} developed an asymptotic theory of portmanteau tests in the standard GARCH framework.
\cite{LKN} suggest a consistent specification test for GARCH$(1,1)$ model. This test is based on a test statistic of Cramér-Von Mises type. Recently, \cite{FWZ-PortTest} proposed a portmanteau test for the Log-GARCH model and the exponential GARCH (EGARCH) model.
\cite{CF-PortTest} work on the APARCH model when the power $\delta$ is known (and thus $\delta$ is not estimated) and suggest a portmanteau test for this class of models.
However, in term of power performance, the authors have showed that: these portmanteau tests are more disappointing since they fail to detect
alternatives of the form $\delta>2$ when the null is $\delta=2$ (see the right array in Table 1 of \cite{CF-PortTest}).
To circumvent the problem, we propose in this work to adopt these portmanteau tests to the case  of APGARCH model when the power $\delta$ is unknown and is jointly estimated with the model's parameters.
%These works do not cover the case where the power is unknown and is estimated with the others parameters.
Consequently, under the null hypothesis of an APGARCH$(p,q)$ model, we shown that the asymptotic distributions of the proposed statistics are a chi-squared distribution as in  \cite{CF-PortTest}.
To obtain this result, we need the following technical (but not restrictive) assumption:

	\indent \textbf{A1:} the support of $\eta_t$ contains at least eleven positive values or eleven negative values.

Notice that \cite{CF-PortTest} need that the support of $\eta_t$ contains at least three positive values or three negative values only.
This is due to the fact that $\delta$ was known in their work.

In Section \ref{estimation}, we recall the results on the quasi-maximum likelihood estimator (QMLE) asymptotic distribution obtained by
\cite{HZ-APGARCH}  when the power $\delta$ is unknown.
Section \ref{Portmanteau test} presents our main aim, which is to complete the work of \cite{CF-PortTest} and to extend the asymptotic theory to the wide class of APGARCH models \eqref{model} when the power $\delta$ is estimated with the other parameters. In Section \ref{Numerical illustration}, we test the null assumption of an APGARCH$(p,q)$ with varying $p$ and $q$ and against different APGARCH models. The null assumption of an APGARCH$(1,1)$ model for different values of $\delta$ are also presented. Section \ref{real datasets} illustrates the portmanteau test for APGARCH models applied to exchange rates. To obtain these results, we use the asymptotic properties obtained by \cite{HZ-APGARCH} for the APGARCH model \eqref{model}.
\section{Quasi-maximum likelihood estimation when the power $\delta$ is unknown}\label{estimation}
Let the parameter space $\Delta \subseteq ]0,+\infty[ \times [0,+\infty[^{2q+p}\times ]0,+\infty[$.

For all $\vartheta=(\omega, \alpha_1^+, \ldots, \alpha_q^+, \alpha_1^-, \ldots, \alpha_q^-, \beta_1,\ldots, \beta_p,\tau)' \in \Delta$, we assume that  $\zeta_t(\vartheta)$ is the strictly stationary and non-anticipative solution of
\begin{equation}\label{zeta}
\zeta_t(\vartheta) = \left(\omega + \sum\limits_{i=1}^q \alpha_i^+(\varepsilon_{t-i}^+)^\tau + \alpha_i^-(-\varepsilon_{t-i}^-)^\tau + \sum\limits_{j=1}^p \beta_j\zeta_{t-j}^\tau(\vartheta)\right)^{1/\tau},
\end{equation}
where $\vartheta$ is equal to an unknown value $\vartheta_0$ belonging to $\Delta$. In the sequel, we let $\zeta_t(\vartheta_0)=\zeta_t$.
Given the realizations $\varepsilon_1,\dots,\varepsilon_n$ (of length $n$) satisfying the APGARCH$(p,q)$ representation \eqref{model},
the variable $\zeta_t(\vartheta)$ can be approximated by $\tilde{\zeta}_t(\vartheta)$ defined recursively by
\begin{equation*}
\tilde{\zeta}_t(\vartheta) =\left( \omega + \sum\limits_{i=1}^q \alpha_{i}^+ (\varepsilon_{t-i}^+)^\tau + \alpha_{i}^- (-\varepsilon_{t-i}^-)^\tau + \sum\limits_{j=1}^p \beta_{j} \tilde{\zeta}_{t-j}^\tau(\vartheta)\right)^{1/\tau},\text{ for }t\geq 1,
\end{equation*}
conditional to the initial values $\varepsilon_0, \ldots, \varepsilon_{1-q}$, $\tilde{\zeta}_0(\vartheta) \geq 0, \ldots, \tilde{\zeta}_{1-p}(\vartheta)\geq 0$.
The quasi-maximum likelihood (QML) method is particularly relevant for GARCH models because
it provides consistent and asymptotically normal estimators for strictly stationary GARCH processes
under mild regularity conditions (but with no moment assumptions on the observed process).
The QMLE is obtained by the standard estimation procedure for GARCH class models. Thus a QMLE of $\vartheta_0$ of the model \eqref{model} is defined as any measurable solution $\hat{\vartheta}_n$ of
\begin{equation}\label{QMLE}
\hat{\vartheta}_n = \underset{\vartheta\in\Delta}{\arg\min} \dfrac{1}{n} \sum\limits_{t=1}^n \tilde{l}_t(\vartheta), \text{ where } \tilde{l}_t(\vartheta) = \dfrac{\varepsilon_t^2}{\tilde{\zeta}_t^2(\vartheta)} + \log(\tilde{\zeta}_t^2(\vartheta)).
\end{equation}
To ensure the asymptotic properties of the QMLE (for the model \eqref{model}) obtained by \cite{HZ-APGARCH}, we need the following assumptions:

	%\indent \textbf{A0:} $(\eta_t)$ is a sequence of iid random variables with $\mathbb{E}\vert \eta_t\vert^r< \infty$ for some $r> 0$.\\
	\indent \textbf{A2:} $\vartheta_0 \in \Delta$ and $\Delta$ is compact.\\
	\indent \textbf{A3:} $\forall \vartheta \in \Delta, \quad \sum_{j=1}^p \beta_j < 1$ and $\gamma(C_{0}) < 0$ where $\gamma(\cdot)$ is the top Lyapunov exponent of the sequence of matrix $C_{0}=\{C_{0t}{,}t\in\mathbb Z\}$ where $C_{0t}$ is defined  in the appendix (see \eqref{cot}). \\
	\indent \textbf{A4:} If $p>0, \mathcal{B}_{\vartheta_0}(z)=1 - \sum_{j=1}^p \beta_{0j}z^j$ has non common root with $\mathcal{A}_{\vartheta_0}^+(z) = \sum_{i=1}^q\alpha_{0i}^+z^i$ and  $\mathcal{A}_{\vartheta_0}^-(z)= \sum_{i=1}^q\alpha_{0i}^-z^i$. Moreover $\mathcal{A}_{\vartheta_0}^+(1) + \mathcal{A}_{\vartheta_0}^-(1) \neq 0$ and $\alpha_{0q}^+ + \alpha_{0q}^- + \beta_{0p} \neq 0$.\\
	\indent \textbf{A5:} $\mathbb{E}[\eta_t^2]=1$ and $\eta_t$ has a positive density on some neighborhood of zero.\\
	\indent \textbf{A6:} $\vartheta_0 \in \overset{\circ}{\Delta}$, where $\overset{\circ}{\Delta}$ denotes the interior of $\Delta$.
	%\indent \textbf{A6:} $\kappa_\eta := \mathbb{E}[\eta_t^4] < \infty$. \\

To ensure the strong consistency of the QMLE, a compactness assumption is required (i.e \textbf{A2}).	
The assumption \textbf{A3} makes reference to the condition of strict stationarity for the model \eqref{model}. Assumptions \textbf{A4} and \textbf{A5}
are made for identifiability reasons and Assumption \textbf{A6} precludes the situation where certain components of $\vartheta_0$ are equal to zero.
%For more details, we invite the reader to see Appendix A in \cite{HZ-APGARCH}.\\	
Then under the assumptions \textbf{A0}, \textbf{A2}--\textbf{A6}, \cite{HZ-APGARCH} showed that $\hat{\vartheta}_n\to\vartheta_0$ $a.s.$ as $n\to\infty$ and
$\sqrt{n}(\hat{\vartheta}_n - \vartheta_0)$ is asymptotically normal with mean 0 and covariance matrix $(\kappa_\eta-1)J^{-1}$, where
$$J := \mathbb{E}_{\vartheta_0}\left[\dfrac{\partial^2 l_t(\vartheta_0)}{\partial \vartheta \partial \vartheta'}\right] =  \mathbb{E}_{\vartheta_0}\left[\dfrac{\partial \log(\zeta_t^2(\vartheta_0))}{\partial \vartheta} \dfrac{\partial \log(\zeta_t^2(\vartheta_0))}{\partial \vartheta'}\right], \text{with }
{l}_t(\vartheta) = \dfrac{\varepsilon_t^2}{{\zeta}_t^2(\vartheta)} + \log({\zeta}_t^2(\vartheta))$$ where $\kappa_\eta := \mathbb{E}[\eta_t^4]< \infty$ by \textbf{A0} and $\zeta_t(\vartheta)$ is given by \eqref{zeta}.
%\begin{propo}{\bf (Hamadeh and Zakoïan, 2011).}\\
%Let $\hat{\vartheta}_n$ a sequence of QMLE satisfying \eqref{QMLE}. Under the above assumptions, we have $\hat{\vartheta}_n\to\vartheta_0$ $a.s.$ as $n\to\infty$. Moreover, it holds
%\begin{equation*}
%\sqrt{n}(\hat{\vartheta}_n - \vartheta_0) \overset{\mathcal{L}}{\longrightarrow} \mathcal{N}(0, (\kappa_\eta-1)J^{-1}),
%\end{equation*}
%where $$J := \mathbb{E}_{\vartheta_0}\left[\dfrac{\partial^2 l_t(\vartheta_0)}{\partial \vartheta \partial \vartheta'}\right] =  4\mathbb{E}_{\vartheta_0}\left[\dfrac{\partial \log(\zeta_t)(\vartheta_0)}{\partial \vartheta} \dfrac{\partial \log(\zeta_t)(\vartheta_0)}{\partial \vartheta'}\right], \text{with }
%{l}_t(\vartheta) = \dfrac{\varepsilon_t^2}{{\zeta}_t^2} + \log({\zeta}_t^2)$$
%and for all $\vartheta \in \Delta$,  $\zeta_t(\vartheta)$ is the strictly stationary and non-anticipative solution of
%\begin{equation}\label{zeta}
%\zeta_t(\vartheta) = \left(\omega + \sum\limits_{i=1}^q \alpha_i^+(\varepsilon_{t-i}^+)^\tau + \alpha_i^-(-\varepsilon_{t-i}^-)^\tau + \sum\limits_{j=1}^p \beta_j\zeta_{t-j}^\tau(\vartheta)\right)^{1/\tau}.
%\end{equation}
%\end{propo}
\section{Portmanteau test}\label{Portmanteau test}
To check the adequacy of a given time series model, for instance an ARMA$(p,q)$ model, it is
common practice to test the significance of the residuals autocorrelations. In the GARCH framework
this approach is not relevant because the process $\eta_t=\varepsilon_t/\zeta_t$ is always a white noise (possibly a
martingale difference) even when the volatility is misspecified.
To check the adequacy of a volatility model, under the null hypothesis
\begin{equation*}
\mathrm{\mathbf{\mathrm{\mathbf{H_0}}}} : \text{the process $(\varepsilon_t)$ satisfies the model \eqref{model}},
\end{equation*}
it is much more fruitful to look at the squared residuals autocovariances
\begin{equation*}
\hat{r}_h = \dfrac{1}{n} \sum\limits_{t = \vert h \vert +1}^n(\hat{\eta}_t^2-1)(\hat{\eta}_{t-\vert h\vert}^2-1),\text{ with } \hat{\eta}_t^2 = \dfrac{\varepsilon_t^2}{\hat{\zeta}_t^2},
\end{equation*}
for $\vert h \vert < n$ and where $\hat{\zeta}_t = \tilde{\zeta}_t(\hat{\vartheta}_n)$ is the quasi-maximum likelihood residuals. For a fixed integer $m\geq 1$, we consider the vector of the first $m$ sample  autocovariances defined by
$$\mathbf{\hat{r}}_m = (\hat{r}_1,\ldots,\hat{r}_m), \text{ such that } 1\leq m < n.$$
Let $I_k$ the identity matrix of size $k$.
The following theorem gives the asymptotic distribution for quadratic forms of autocovariances
of squared residuals.
\begin{thm}\label{resultat.r_h}
Under the assumptions \textbf{A0}--\textbf{A6}, if $(\varepsilon_t)$ is the non-anticipative
and stationary solution of the APGARCH$(p,q)$ model \eqref{model}, then, when $n\to\infty$, we have
\begin{equation*}
\sqrt{n} \mathbf{\hat{r}}_m \overset{\mathcal{L}}{\longrightarrow} \mathcal{N}\left(0,D\right)\text{ where }
{D}=({\kappa}_\eta - 1)^2 I_m - ({\kappa}_\eta-1){C}_m{J}^{-1}{C}'_m
\end{equation*}
is nonsingular and where the matrix ${C}_m$ is given by \eqref{C_m} in the proof of Theorem~\ref{resultat.r_h}.
\end{thm}
The proof of this result is postponed to Section \ref{proof}.

The standard portmanteau test for checking that the data is a realization of a strong white noise is that \cite{BP} or \cite{LB}.
Both of these tests are based on the residuals autocorrelations $\hat\rho(h)$  and they are defined by
\begin{equation}\label{bp}
 Q^{\textsc{bp}}_m=n\sum_{h=1}^m\hat\rho^2(h) \text{ and }   {Q}_m^{\textsc{lb}}=n(n+2)\sum_{h=1}^m\frac{\hat\rho^2(h)}{n-h},
\end{equation}
where $n$ is the length of the series and $m$ is a fixed integer.
Under the assumption that the  noise sequence is iid, the standard test procedure consists in rejecting the strong white noise hypothesis
if the statistics (\ref{bp})  are larger than a certain quantile of a chi-squared distribution.
These tests are not robust to conditional heteroscedasticity  or other processes displaying
a second order dependence. Indeed such nonlinearities
may arise for instance when the observed process $(\varepsilon_t)$ follows a GARCH representation.
 Other situations where the standard  tests are not robust can be found for instance in  \cite{frz} or \cite{yac},
 who showed that: for an ARMA model with uncorrelated but dependent noise process, the asymptotic distributions of
 the statistics defined in (\ref{bp}) are no longer chi-squared distributions but a mixture of chi-squared distributions.
In the APGARCH framework, we may wish to simultaneously test the nullity of the first $m$ autocovariances using more robust
portmanteau statistics.

In order to state our second result, we also need further notations. Let $\hat{\kappa}_\eta$, $\hat{J}$ and $\hat{C}_m$ be weakly consistent estimators of
$\kappa_\eta$, $J$ and ${C}_m$ involved in the asymptotic normality of $\sqrt{n} \mathbf{\hat{r}}_m$ (see Theorem~\ref{resultat.r_h}).
For instance, $\kappa_\eta$ and $J$ can be estimated by their empirical or observable counterparts
given by
\begin{equation*}
\hat{\kappa}_\eta = \dfrac{1}{n} \sum\limits_{t=1}^n \dfrac{\varepsilon_t^4}{\tilde{\zeta}_t^4(\hat{\vartheta}_n)} \qquad \mbox{and}\qquad \hat{J} = \dfrac{1}{n}\sum\limits_{t=1}^n \dfrac{\partial \log\tilde{\zeta}_t^2(\hat{\vartheta}_n)}{\partial \vartheta}\dfrac{\partial \log\tilde{\zeta}_t^2(\hat{\vartheta}_n)}{\partial \vartheta'}.
\end{equation*}
We can write the vector of parameters $\vartheta := (\theta', \tau)'$ where $\theta\in\mathbb{R}^{2q+p+1}$ depends on the coefficients $\omega, \alpha_1^+, \ldots, \alpha_q^+, \alpha_1^-, \ldots, \alpha_q^-, \beta_1,\ldots, \beta_p$. The parameter of interest becomes $\vartheta_0 :=(\theta_0',\delta)'$, where $\theta_0= (\omega_0, \alpha_{01}^+, \ldots, \alpha_{0q}^+, \alpha_{01}^-, \ldots, \alpha_{0q}^-, \beta_{01},\ldots, \beta_{0p})'$. With the previous notation, for all $\vartheta= (\theta', \tau)'\in\Delta$, the derivatives in the expression of $\hat{J}$ can be recursively computed for $t> 0$ by
\begin{equation*}
\begin{aligned}
\dfrac{\partial \tilde{\zeta}_t^\tau(\vartheta)}{\partial \theta} &= \underline{\tilde c}_t(\vartheta) + \sum\limits_{j=1}^p \beta_j \dfrac{\partial \tilde{\zeta}_{t-j}^\tau(\vartheta)}{\partial \theta},\\
\dfrac{\partial \tilde{\zeta}_t^\tau(\vartheta)}{\partial \tau} &= \sum\limits_{i=1}^q \alpha_i^+\log(\varepsilon_{t-i}^+)(\varepsilon_{t-i}^+)^\tau + \alpha_i^-\log(-\varepsilon_{t-i}^-)(-\varepsilon_{t-i}^-)^\tau  + \sum\limits_{j=1}^p \beta_j \dfrac{\partial \tilde{\zeta}_{t-j}^\tau(\vartheta)}{\partial \tau},
\end{aligned}
\end{equation*}
with the initial values $\partial \tilde{\zeta}_t(\vartheta)/\partial \vartheta = 0$, for all $t = 0, \ldots, 1-p$ and
\begin{equation}\label{ctilde}
\underline{\tilde c}_t(\vartheta) = (1, (\varepsilon_{t-1}^+)^\tau, \ldots, (\varepsilon_{t-q}^+)^\tau, (-\varepsilon_{t-1}^-)^\tau, \ldots, (-\varepsilon_{t-q}^-)^\tau, \tilde{\zeta}_{t-1}^\tau(\vartheta), \ldots, \tilde{\zeta}_{t-p}^\tau(\vartheta))'.
\end{equation}
By convention, $\log(\varepsilon_{t}^+) = 0$ if $\varepsilon_t \leq 0$ and respectively $\log(-\varepsilon_{t}^-) = 0$ if $\varepsilon_t \geq 0$.

For the matrix $\hat{C}_m$ of size $m \times (2q+p+2)$, one can take
\begin{equation}
\hat{C}_m(h,k) = -\dfrac{1}{n}\sum\limits_{t=h+1}^n(\hat{\eta}^2_{t- h} -1) \dfrac{1}{\tilde{\zeta}_t^2(\hat{\vartheta}_n)} \dfrac{\partial \tilde{\zeta}_t^2(\hat{\vartheta}_n)}{\partial \vartheta_k} \text{ for }1\leq h \leq m\text{ and }1\leq k \leq 2q+p+2,
\end{equation}
where $\hat{C}_m(h,k)$ denotes the $(h,k)$ element of the matrix $\hat{C}_m$. Let
$\hat{D} =(\hat{\kappa}_\eta-1)^2 I_m-(\hat{\kappa}_\eta-1)\hat{C}_m\hat{J}^{-1}\hat{C}_m$ be a weakly consistent estimator of the matrix $D$.
The following result is established in the case where the power is unknown and estimated with the others parameters.
\begin{thm}\label{resultat}
Under Assumptions of Theorem~\ref{resultat.r_h} and $\mathrm{\mathbf{H_0}}$, when $n\to\infty$, we have
\begin{equation*}
n \mathbf{\hat{r}}_m'\hat{D}^{-1}\mathbf{\hat{r}}_m \overset{\mathcal{L}}{\longrightarrow} \chi^2_m.
\end{equation*}
\end{thm}
The proof of this result is postponed to Section \ref{proof}.

The adequacy of the APGARCH$(p,q)$ model \eqref{model} is then rejected at the asymptotic level $\alpha$ when
\begin{equation}
n \mathbf{\hat{r}}_m'\hat{D}^{-1}\mathbf{\hat{r}}_m\  >\  \chi^2_m(1-\alpha),
\end{equation}
where $\chi^2_m(1-\alpha)$ represents the $(1-\alpha)-$quantile of the chi-square distribution with $m$ degrees of freedom.
\section{Numerical illustration}\label{Numerical illustration}
By means of Monte Carlo experiments, we investigate the finite sample properties of the test introduced in this paper.
The numerical illustrations of this section are made  with the free
statistical software RStudio  (see https://www.rstudio.com) in Rcpp language.
We simulated $N=1,000$ independent replications of size $n=500$ and $n=5,000$ of the APGARCH$(p,q)$ model \eqref{model} with the orders $(p,q)\in\{0,1,2\}\times\{1,2,3\}$. The distribution of $\eta_t$ is a Student law with $9$ degrees of freedom, standardized to obtain a variance equal to 1.
%The vector of parameters used to simulate an APGARCH$(1,1)$ is $\vartheta_{0,11} := (0.04, 0.02, 0.13, 0.85, 1)'$. For the other models the parameters $\vartheta_{0,pq}$ are close to $\vartheta_{0,11}$ in such a way that the existence of fourth order moment is satisfied.
%The size of the observations is $n=5,000$, which could correspond in practice to the length for daily financial series or higher-frequency data.

For each of these $N$ replications and each APGARCH$(p,q)$ models considered,  we use the QMLE method to estimate the corresponding coefficients $\vartheta_{0,pq}$ and we apply portmanteau test to the squared residuals for different values of $m$,  where $m$ is the number of autocorrelations used in the portmanteau test statistic.
 At the nominal level $\alpha=5\%$ the confidence interval of the nominal level is $[3.6\%, 6.4\%]$ with a probability $95\%$ and $[3.2\%,6.8\%]$ with a probability $99\%$.

The left array in Tables 1 (resp. Tables 2) represents the number of rejection in percentage of the orders $p$ and $q$ for the corresponding  APGARCH$(p,q)$ models
for $n=500$ (resp. $n=5,000$). These tests are done for the nominal level $\alpha=5\%$. As excepted, for all models the percentages of rejection belongs to the confident interval with  probabilities $95\%$ and $99\%$, except for  the APGARCH$(2,1)$ and APGARCH$(1,3)$ models when $=500$ and $m\geq4$. Consequently the proposed test well controls the error of first kind for the candidates models
when the number of observations is $n=5,000$, which could correspond in practice to the length for daily financial series or higher-frequency data.

We now study the empirical power under the null of an APGARCH$(1,2)$ model.
The right array in Tables 1 (resp. Tables 2) displays  the relative rejection frequencies (also in percentage) over the $N$ independent replications in the case that the null is an APGARCH$(1,2)$ model
for $n=500$ (resp. $n=5,000$). In these cases, we also estimate the power $\tau$ of different models with true value $\delta=1$, which correspond to the TGARCH models of \cite{RZ-TGARCH}. The test makes the difference between the models when the size $n$ increases (see right  array in Table 2).%, but it seems to confuse the model with the APGARCH$(1,2)$ (see right  array in Table 1). %In others simulations, we reproduce the right array in table 1 under the null assumption that the model is an APGARCH$(1,1)$ and we can see that the test also confuse the model APGARCH$(1,1)$ and APGARCH$(2,1)$.\\

\cite{CF-PortTest} work on the APARCH model when the power $\delta$ is known and suggest a portmanteau test for this class of models.
However, in term of power performance, the authors have showed that: these portmanteau tests are more disappointing since they fail to detect
alternatives of the form $\delta>2$ when the null is $\delta=2$ (see the right array in Table 1 of \cite{CF-PortTest}).
Contrary to \cite{CF-PortTest}, we estimate the power $\delta$ and consequently we can not compare our simulations.
Nevertheless, in Table 3 we present the frequencies of rejection in percentage for the model APGARCH$(1,1)$ when the power $\delta\in [0.5{,}3]$ is estimated.
To simulate the different trajectories, we use the  parameter $\theta_0 = (0.04, 0.02, 0.13, 0.85)'$ used by \cite{CF-PortTest}.
%For all the power tested, the results are the most of the time in the confident interval.
However, from Table 3 the test do not reject the null hypothesis when $\delta$ is higher than $2$.
So, this problem seems to be overcome when the power $\delta$ is unknown and is jointly estimated with the model's parameters.
We draw the same conclusion that the test also controls well the error of first kind at different asymptotic level $\alpha$.

\begin{center}
\begin{tabular}{p{1.3cm}|m{0.7cm}cccccc||cccccc}
\multicolumn{8}{c}{Empirical Size} & \multicolumn{6}{c}{Empirical Power}\\
\hline \hline
\centering level & \multirow{2}{*}{\centering$(p,q)$} & \multicolumn{6}{c||}{$m$} & \multicolumn{6}{c}{$m$} \\
\cline{3-14}
&& 2 & 4 & 6 & 8 & 10 & 12 & 2 & 4 & 6 & 8 & 10 & 12\\
\hline
&\centering$(0,1)$ & 4.8 & 5.7 & 4.5 & 6.3 & 5.7 & 4.6 & 15.4 & 19.0 & 18.4 & 15.7 & 16.6 & 13.6 \\
&\centering$(1,1)$ & 3.7 & 4.3 & 5.1 & 6.0 & 5.7 & 6.8 & 6.6 & 7.7 & 7.3 & 8.8 & 8.2 & 7.8 \\
&\centering$(1,2)$ & 4.3 & 4.2 & 6.1 & 4.7 & 4.1 & 4.3 & 4.3 & 4.2 & 6.1 & 4.7 & 4.1 & 4.3  \\
&\centering$(1,3)$ & 5.8 & 7.1 & 10.1 & 8.1 & 9.6 & 8.9 & 10.6 & 13.2 & 15.2 & 10.6 & 12.5 & 12.2\\
&\centering$(2,1)$ & 6.8 & 8.0 & 8.5 & 10.3 & 8.9 & 8.9 & 10.1 & 9.8 & 9.6 & 10.9 & 9.2 & 8.7 \\
\hline
\end{tabular}
\end{center}
\begin{center}
\footnotesize{Table 1: Relative frequencies (in $\%$) of rejection when $n=500$\\ \vspace{-3mm}
\begin{flushleft}
Left part: Relative frequencies of rejection for different APGARCH$(p,q)$ models with the power estimated. \\
Right part: Relative frequencies of rejection when the model is an APGARCH$(1,2)$ with $\vartheta_{0,12} = (0.04, 0.02,0.005, 0.13, 0.05, 0.6, 1)'$.
\end{flushleft}}\end{center}
\vspace*{1mm}

\begin{center}
\begin{tabular}{p{1.3cm}|m{0.7cm}cccccc||cccccc}
\multicolumn{8}{c}{Empirical Size} & \multicolumn{6}{c}{Empirical Power}\\
\hline \hline
\centering level & \multirow{2}{*}{\centering$(p,q)$} & \multicolumn{6}{c||}{$m$} & \multicolumn{6}{c}{$m$} \\
\cline{3-14}
&& 2 & 4 & 6 & 8 & 10 & 12 & 2 & 4 & 6 & 8 & 10 & 12\\
\hline
&\centering$(0,1)$ & 4.5 & 4.7 & 4.2 & 5.9 & 6.3 & 5.5 & 99.6 & 99.9 & 99.4 & 99.7 & 99.6 & 99.0 \\
&\centering$(1,1)$ & 4.6 & 6.5 & 4.6 & 4.8 & 7.2 & 6.3 & 24.2 & 18.6 & 16.4 & 14.6 & 11.6 & 9.8 \\
&\centering$(1,2)$ & 4.4 & 5.5 & 5.7 & 4.8 & 5.2 & 5.2 & 4.4 & 5.5 & 5.7 & 4.8 & 5.2 & 5.2 \\
&\centering$(1,3)$ & 5.5 & 6.5 & 6.2 & 6.3 & 6.2 & 8.0 & 10.6 & 14.5 & 14.5 & 14.4 & 13.4 & 11.9 \\
&\centering$(2,1)$ & 4.7 & 6.2 & 4.7 & 6.3 & 6.1 & 6.5 & 42.2 & 38.7 & 35.2 & 33.0 & 32.5 & 28.5 \\
\hline
\end{tabular}
\end{center}
\begin{center}
\footnotesize{Table 2: Relative frequencies (in $\%$) of rejection when $n=5,000$\\ \vspace{-3mm}
\begin{flushleft}
Left part: Relative frequencies of rejection for different APGARCH$(p,q)$ models with the power estimated. \\
Right part: Relative frequencies of rejection when the model is an APGARCH$(1,2)$, with $\vartheta_{0,12} = (0.04, 0.02,0.005, 0.13, 0.05, 0.6, 1)'$.
\end{flushleft}}\end{center}
\vspace*{1mm}

%
%
%
%\  \\
%\hspace{-0.7cm}
\begin{minipage}[c]{8.7cm}
      \begin{tabular}{p{1.3cm}|m{0.7cm}cccccc}
      	\hline \hline
     	\centering level &  \multirow{2}{*}{\centering$\delta$} & \multicolumn{5}{c}{$m$}\\
      	\cline{3-8}
      	&& 2 & 4 & 6 & 8 & 10 & 12\\
	\hline
	\multirow{6}{*}{$\alpha = 1\%$} & \centering0.5 & 1.0 & 1.2 & 1.1 & 1.2 & 2.3 & 2.4 \\
	&\centering1 & 1.3 & 1.3 & 1.4 & 1.7 & 1.3 & 1.7 \\
	&\centering1.5 & 1.4 & 1.2 & 1.0 & 1.7 & 1.3 & 1.7 \\
	&\centering2 & 1.8 & 1.4 & 1.0 & 1.4 & 1.8 & 1.7 \\
	&\centering2.5 & 1.1 & 1.1 & 0.9 & 1.8 & 1.8 & 1.1 \\
	&\centering3 & 1.3 & 1.9 & 1.2 & 1.7 & 1.7 & 1.6 \\

	\hline
      \end{tabular}
\end{minipage}
\begin{minipage}[c]{9cm}
      \begin{tabular}{p{1.3cm}|m{0.7cm}cccccc}
      	\hline \hline
     	\centering level & \multirow{2}{*}{\centering$\delta$} & \multicolumn{5}{c}{$m$}\\
      	\cline{3-8}
      	&& 2 & 4 & 6 & 8 & 10 & 12\\
	\hline
	\multirow{6}{*}{${\alpha = 5\%}$} & \centering0.5 & 4.7 & 5.6 & 4.8 & 6.0 & 6.0 & 6.4 \\
	&\centering1 & 3.2 & 4.4 & 5.8 & 5.4 & 5.4 & 5.1 \\
	&\centering1.5 & 4.4 & 4.2 & 5.5 & 6.0 & 4.6 & 5.1 \\
	&\centering2 & 5.2 & 4.7 & 5.5 & 5.1 & 5.7 & 7.4 \\
	&\centering2.5 & 3.7 & 4.9 & 4.9 & 4.3 & 5.0 & 5.3 \\
	&\centering3 & 3.8 & 3.4 & 4.8 & 5.4 & 4.9 & 6.6\\
	\hline
      \end{tabular}
\end{minipage}
\begin{center}
\footnotesize{
Table 3: \\ \vspace{-3mm}
\begin{flushleft}
Relative frequencies (in $\%$) of rejection for an APGARCH$(1,1)$ model with different power coefficients
and $\vartheta_{0,11} := (0.04, 0.02, 0.13, 0.85, \delta)'$.
Left part: the nominal level is $\alpha = 1\%$ and $\alpha = 5\%$ in the right part.
\end{flushleft}}
\end{center}
\section{Adequacy of APGARCH models for real datasets}\label{real datasets}
%By practitioners, the GARCH$(1,1)$ model is the most widely used to estimate
%the volatility of daily returns. In general, this model is chosen a priori, without implementing
%any statistical identification procedure. This practice is motivated by the common belief that the
%GARCH$(1,1)$ (or its simplest asymmetric extensions) is sufficient to capture the properties of
%daily financial series and that higher-order models may be unnecessarily complicated.
%We will show that, for a large number of series, this practice is not always statistically justified.
We consider the daily return of four exchange rates EUR/USD (Euros Dollar), EUR/JPY (Euros Yen), EUR/GBP (Euros Pounds) and EUR/CAD (Euros Canadian dollar). The observations covered the period from November 01, 1999 to April 28, 2017 which correspond to $n=4,478$ observations. The data were obtain from the website of the National Bank of Belgium (https://www.nbb.be).

Table 4 displays the $p-$values for adequacy of the APGARCH$(p,q)$  for daily returns of exchange rates  based on $m$ squared residuals autocovariances, as well as the estimated power.
The APGARCH$(0,1)$ model assumption is rejected for each series and is not adapted to these kinds of series.
The APGARCH$(1,2)$ model is rejected for EUR/GBP and EUR/CAD whereas the APGARCH$(1,1)$ and APGARCH$(2,1)$ models seem the most appropriate for the exchange rates. The APGARCH$(2,2)$ model assumption is only rejected for the exchange rates EUR/CAD.
From the last column of Table 4, we can also see that the estimated power $\hat\tau$ is not necessary equal to 1 or 2 and is different for each series.

%The AIC criterion tends to always select the APGARCH$(1,1)$.
The portmanteau test is thus an important tool in the validation process.
From the empirical results and the simulation experiments,
we draw the conclusion that the proposed portmanteau test based on squared
residuals of an APGARCH$(p,q)$ (when the power is unknown and is jointly estimated with the model's parameters) is efficient to detect a misspecification of the order $(p,q)$.
% but the test appears to have trouble to distinguish some sub-class of models as the APGARCH$(1,1)$ and the APGARCH$(2,1)$ or the APGARCH$(2,1)$ and the APGARCH$(2,2)$.

\begin{center}
\begin{tabular}{p{0.7cm}c@{\hskip0.25cm}c@{\hskip0.25cm}c@{\hskip0.25cm}c@{\hskip0.25cm}c@{\hskip0.25cm}c@{\hskip0.25cm}c@{\hskip0.25cm}c@{\hskip0.25cm}c@{\hskip0.25cm}c@{\hskip0.25cm}c@{\hskip0.25cm}c@{\hskip0.25cm}c}
\hline\hline
& \multicolumn{12}{c}{$m$} & \multirow{2}{*}{$\hat{\tau}$}\\
\cline{2-13}
& 1 & 2 & 3 & 4 & 5 & 6 & 7 & 8 & 9 & 10 & 11 & 12 &\\
\hline
\\
\multicolumn{13}{l}{Portmanteau tests for adequacy of the APGARCH(0,1)}\\
USD & 0.009 & 0.003 & 0.000 & 0.000 & 0.000 & 0.000 & 0.000 & 0.000 & 0.000 & 0.000 & 0.000 & 0.000 & 1.77 \\
JPY & 0.160 & 0.000 & 0.000 & 0.000 & 0.000 & 0.000 & 0.000 & 0.000 & 0.000 & 0.000 & 0.000 & 0.000 & 1.23 \\
GBP & 0.697 & 0.000 & 0.000 & 0.000 & 0.000 & 0.000 & 0.000 & 0.000 & 0.000 & 0.000 & 0.000 & 0.000 & 1.98 \\
CAD & 0.049 & 0.000 & 0.000 & 0.000 & 0.000 & 0.000 & 0.000 & 0.000 & 0.000 & 0.000 & 0.000 & 0.000 & 2.35 \\
\\
\multicolumn{13}{l}{Portmanteau tests for adequacy of the APGARCH(1,1)}\\
USD & 0.888 & 0.533 & 0.715 & 0.671 & 0.764 & 0.814 & 0.687 & 0.704 & 0.788 & 0.817 & 0.874 & 0.906 & 1.05 \\
JPY & 0.113 & 0.261 & 0.442 & 0.605 & 0.735 & 0.442 & 0.550 & 0.578 & 0.591 & 0.342 & 0.401 & 0.478 & 1.11\\
GBP & 0.037 & 0.087 & 0.181 & 0.166 & 0.242 & 0.346 & 0.362 & 0.292 & 0.377 & 0.410 & 0.406 & 0.490 & 1.33 \\
CAD & 0.027 & 0.078 & 0.157 & 0.254 & 0.174 & 0.254 & 0.291 & 0.269 & 0.346 & 0.435 & 0.517 & 0.536 & 1.56 \\
\\
\multicolumn{13}{l}{Portmanteau tests for adequacy of the APGARCH(1,2)}\\
USD & 0.673 & 0.489 & 0.672 & 0.648 & 0.739 & 0.780 & 0.647 & 0.646 & 0.739 & 0.767 & 0.832 & 0.870 & 1.08 \\
JPY & 0.003 & 0.009 & 0.025 & 0.051 & 0.089 & 0.055 & 0.086 & 0.098 & 0.121 & 0.063 & 0.083 & 0.114 & 1.11\\
GBP & 0.000 & 0.000 & 0.001 & 0.001 & 0.003 & 0.006 & 0.007 & 0.005 & 0.009 & 0.011 & 0.015 & 0.023 & 1.33 \\
CAD & 0.000 & 0.000 & 0.000 & 0.000 & 0.000 & 0.000 & 0.000 & 0.000 & 0.000 & 0.000 & 0.000 & 0.000 & 1.53 \\
\\
\multicolumn{13}{l}{Portmanteau tests for adequacy of the APGARCH(2,1)}\\
USD & 0.471 & 0.544 & 0.682 & 0.622 & 0.733 & 0.787 & 0.651 & 0.659 & 0.750 & 0.781 & 0.843 & 0.877 & 1.05 \\
JPY & 0.379 & 0.680 & 0.855 & 0.941 & 0.977 & 0.657 & 0.763 & 0.796 & 0.747 & 0.342 & 0.294 & 0.351 & 1.10 \\
GBP & 0.193 & 0.362 & 0.566 & 0.455 & 0.564 & 0.687 & 0.689 & 0.587 & 0.676 & 0.696 & 0.669 & 0.746 & 1.34 \\
CAD & 0.170 & 0.277 & 0.440 & 0.594 & 0.403 & 0.523 & 0.567 & 0.515 & 0.607 & 0.698 & 0.768 & 0.779 & 1.61 \\
\\
\multicolumn{13}{l}{Portmanteau tests for adequacy of the APGARCH(2,2)}\\
USD & 0.849 & 0.448 & 0.630 & 0.600 & 0.715 & 0.784 & 0.634 & 0.693 & 0.779 & 0.815 & 0.870 & 0.907 & 1.02 \\
JPY & 0.057 & 0.154 & 0.291 & 0.439 & 0.579 & 0.285 & 0.387 & 0.437 & 0.434 & 0.217 & 0.243 & 0.304 & 1.10 \\
GBP & 0.008 & 0.016 & 0.034 & 0.033 & 0.050 & 0.081 & 0.107 & 0.095 & 0.136 & 0.166 & 0.167 & 0.220 & 1.34 \\
CAD & 0.000 & 0.000 & 0.000 & 0.000 & 0.000 & 0.000 & 0.000 & 0.000 & 0.000 & 0.000 & 0.000 & 0.000 & 1.61\\
\hline
\end{tabular}
\end{center}
\begin{center}
\footnotesize{Table 4: Portmanteau test $p-$values for adequacy of the APGARCH$(p,q)$  for daily returns of exchange rates, based on $m$ squared residuals autocovariances.}
\end{center}

\section{Appendix : Proofs}
\label{proof}
We recall that for all $\vartheta \in \Delta$, $\zeta_t(\vartheta)$ is the strictly stationary and non-anticipative solution of \eqref{zeta}.

The matrix $J$ can be rewritten as
\begin{equation*}
J = \mathbb{E}_{\vartheta_0}\left[ \dfrac{1}{\zeta_{t}^4(\vartheta_0)} \dfrac{\partial \zeta_{t}^2(\vartheta_0)}{\partial \vartheta} \dfrac{\partial \zeta_{t}^2(\vartheta_0)}{\partial \vartheta'}\right].
\end{equation*}
First, we shall need some technical results which are essentially contained in  \cite{HZ-APGARCH}.
Let $K$ and $\rho$ be generic constants, whose values will be modified along the proofs, such that $K>0$ and $\rho\in]0,1[$.
\subsection{Reminder on technical issues on quasi likelihood method for APGARCH models}
The starting point is the asymptotic irrelevance of the initial values. Under  \textbf{A0}, \textbf{A2}--\textbf{A6},
\cite{HZ-APGARCH} show that:
\begin{equation}\label{res0}
\sup\limits_{\vartheta\in\Delta} | \zeta_{t}^\tau(\vartheta) - \tilde{\zeta}_{t}^\tau(\vartheta) | \leq K\rho^t.
\end{equation}
Similar properties also hold for the derivatives with respect to $\vartheta$ of $\zeta_{t}^\tau(\vartheta) - \tilde{\zeta}_{t}^\tau(\vartheta)$.
We sum up the properties that we shall need in the sequel. We refer to  \cite{HZ-APGARCH}  for a more detailed treatment.
For some $s \in ]0,1[$, we have
\begin{equation}\label{res1}
\mathbb{E}\vert \varepsilon_0 \vert^{2s} < \infty, \qquad \mathbb{E} \sup\limits_{\vartheta\in\Delta} | \zeta_{t}^{2s} | < \infty, \qquad \mathbb{E} \sup\limits_{\vartheta\in\Delta} | \tilde{\zeta}_{t}^{2s} | < \infty.
\end{equation}
Moreover, from \eqref{res0}, the mean-value theorem implies that
\begin{equation}\label{res2}
\sup\limits_{\vartheta\in\Delta} | \zeta_{t}^2(\vartheta) - \tilde{\zeta}_{t}^2(\vartheta)| \leq K\rho^t\sup\limits_{\vartheta\in\Delta}\max\{\zeta_t^2(\vartheta), \tilde{\zeta}_t^2(\vartheta)\}.
\end{equation}
For all $d\geq 1$
\begin{equation}\label{res3}
\mathbb{E}\left\Vert \sup\limits_{\vartheta\in\Delta} \dfrac{1}{\zeta_{t}^\tau(\vartheta)}\dfrac{\partial \zeta_t^\tau(\vartheta)}{\partial \vartheta}\right\Vert^d <\infty, \qquad \mathbb{E}\left\Vert \sup\limits_{\vartheta\in\Delta} \dfrac{1}{\zeta_{t}^\tau(\vartheta)}\dfrac{\partial^2 \zeta_t^\tau(\vartheta)}{\partial \vartheta \partial \vartheta'}\right\Vert^d <\infty.
\end{equation}
There exists a neighborhood $\mathcal{V}(\vartheta_0)$ of $\vartheta_0$ such that for all $\xi > 0$ and $a = 1-(\delta/\tau)(1-s)> 0$
\begin{equation*}
\sup\limits_{\vartheta\in\mathcal{V}(\vartheta_0)} \left(\dfrac{\zeta_t^2(\vartheta_0)}{\zeta_t^2(\vartheta)}\right) \leq \left(K + K\sum\limits_{i=1}^q\sum\limits_{k=0}^\infty (1 + \xi)^k\rho^{ak}\vert \varepsilon_{t-i-k}\vert^{2\tau}\right)^{2/\tau},
\end{equation*}
and it holds that
\begin{equation}\label{res4}
\mathbb{E}\left\vert \sup\limits_{\vartheta\in\mathcal{V}(\vartheta_0)} \left(\dfrac{\zeta_t^2(\vartheta_0)}{\zeta_t^2(\vartheta)}\right)\right\vert < \infty.
\end{equation}
The matrix $J$ is invertible and
\begin{equation}\label{res5}
\quad  \sqrt{n}(\hat{\vartheta}_n - \vartheta_0) = J^{-1} \dfrac{1}{\sqrt{n}} \sum\limits_{t=1}^ns_t \dfrac{1}{\zeta_t^2}\dfrac{\partial \zeta_t^2(\vartheta_0)}{\partial \vartheta} + \mathrm{o}_{\mathbb P}(1),\quad \text{with } s_t = \eta_t^2-1.
\end{equation}
\subsection{Proof of Theorem~\ref{resultat.r_h}}
The proof of Theorem \ref{resultat.r_h} is close to the proof of \cite{CF-PortTest}. Only the invertibility of the matrix $D$ needs to be adapted.
But, to understand the proofs and to have its own autonomy, we rewrite all the proof. We also decompose this proof in 3 following steps.
\begin{enumerate}[$\qquad (i)$]
	\item Asymptotic impact of the unknown initial values on the statistic $\mathbf{\hat{r}}_m$.
	\item Asymptotic distribution of $\sqrt{n}\mathbf{\hat{r}}_m$.
	\item Invertibility of the matrix $D$.
	%\item Convergence to the matrix $D$.
\end{enumerate}
We now introduce the vector of $m$ autocovariances $\mathbf{r}_m = (r_1, \ldots, r_m)'$ where the $h$-th element is define as
\begin{equation*}
r_h = \dfrac{1}{n}\sum\limits_{t=h+1}^n s_ts_{t-h}\ , \quad \text{with } s_t = \eta_t^2-1 \text{ and } 0<h < n.
\end{equation*}
Let $s_t(\vartheta)=\eta^2_t(\vartheta)-1$ with $\eta_t(\vartheta) = \varepsilon_t / \zeta_t(\vartheta)$ and $\tilde{s}_t(\vartheta)= \tilde\eta^2_t(\vartheta)-1$ with $\tilde\eta_t(\vartheta) = \varepsilon_t / \tilde\zeta_t(\vartheta)$.
Let $r_h(\vartheta)$ obtained by replacing $\eta_t$ by $\eta_t(\vartheta)$ in $r_h$ and $\tilde{r}_h(\vartheta)$ by replacing ${\eta}_t$ by $\tilde{\eta}_t(\vartheta)$ in $r_h$. The vectors $\mathbf{r}_m(\vartheta) = (r_1(\vartheta),\ldots, r_m(\vartheta))'$ and $\tilde{\mathbf{r}}_m(\vartheta) = (\tilde{r}_1(\vartheta),\ldots, \tilde{r}_m(\vartheta))'$ are such that $\mathbf{r}_m = \mathbf{r}_m(\theta_0)$, $\tilde{\mathbf{r}}_m = \tilde{\mathbf{r}}_m(\theta_0)$ and $\hat{\mathbf{r}}_m = \tilde{\mathbf{r}}_m(\hat{\vartheta}_n)$.

\textbf{$(i)$  Asymptotic impact of the unknown initial values on the statistic $\mathbf{\hat{r}}_m$.}

We have $s_t(\vartheta)s_{t-h}(\vartheta) - \tilde{s}_t(\vartheta)\tilde{s}_{t-h}(\vartheta) = a_t + b_t$ with $a_t = \{s_t(\vartheta) - \tilde{s}_t(\vartheta)\}s_{t-h}(\vartheta)$ and $b_t = \tilde{s}_t(\vartheta)\{s_{t-h}(\vartheta) - \tilde{s}_{t-h}(\vartheta)\}$.
Using \eqref{res2} and $\inf_{\vartheta \in \Delta} \tilde{\zeta}^2_t \geq \inf_{\vartheta \in \Delta}\omega^{2/\tau} > 0$, we have
\begin{equation*}
\vert a_t \vert+\vert b_t \vert \leq K \rho^t \varepsilon_t^2(\varepsilon_{t-h}^2 + 1) \sup\limits_{\vartheta \in \Delta} \max\{ \tilde{\zeta}_t^2, \zeta_t^2\} \ .
\end{equation*}
Using the inequality $(a+b)^s \leq a^s + b^s$, for $a,b \geq 0$ and $s\in ]0,1[$, \eqref{res1} and Hölder's inequality, %and the inequality $(a+b)^s \leq a^s + b^s$, for $a,b \geq 0$ and $s\in ]0,1[$,
we have for some $s^\ast \in ]0,1[$ sufficiently small
\begin{equation*}
\mathbb{E}\left\vert \dfrac{1}{\sqrt{n}}\sum\limits_{t=1}^n \sup\limits_{\vartheta\in\Delta}\vert a_t \vert \right\vert^{s^\ast} \leq K \dfrac{1}{n^{s^\ast/2}}\sum\limits_{t=1}^n \rho^{ts^\ast} \underset{n\to \infty}{\longrightarrow} 0.
\end{equation*}
We deduce that $n^{-1/2}\sum_{t=1}^n \sup_{\vartheta \in \Delta} \vert a_t \vert = \mathrm{o}_{\mathbb P}(1)$. We have the same convergence for $b_t$, and for the derivatives of $a_t$ and $b_t$.
Consequently, we obtain
\begin{equation}\label{11}
\sqrt{n}\Vert \mathbf{r}_m - \mathbf{\tilde{r}}_m \Vert = \mathrm{o}_{\mathbb P}(1), \qquad \sup\limits_{\vartheta\in\Delta}\left\Vert \dfrac{\partial \mathbf{r}_m}{\partial\vartheta} - \dfrac{\partial \mathbf{\tilde{r}}_m}{\partial\vartheta}\right\Vert = \mathrm{o}_{\mathbb P}(1),\text{ as }n\to\infty.
\end{equation}
The unknown initial values have no asymptotic impact on the statistic $\mathbf{\hat{r}}_m$.

\textbf{$(ii)$ Asymptotic distribution of $\sqrt{n}\mathbf{\hat{r}}_m$.}

We now show that the asymptotic distribution of $\sqrt{n}\mathbf{\hat{r}}_m$ is deduced from the joint distribution of $\sqrt{n}\mathbf{r}_m$ and of the QMLE.
Using \eqref{11} and a Taylor expansion of $\mathbf{r}_m(\cdot)$ around $\hat{\vartheta}_n$ and $\vartheta_0$, we obtain
\begin{equation*}
\begin{aligned}
\sqrt{n} \mathbf{\hat{r}}_m &= \sqrt{n}\mathbf{\tilde{r}}_m(\vartheta_0) + \dfrac{\partial \mathbf{\tilde{r}}_m(\vartheta^\ast)}{\partial \vartheta} \sqrt{n}(\hat{\vartheta}_n - \vartheta_0) \\
&= \sqrt{n}\mathbf{r}_m + \dfrac{\partial \mathbf{r}_m(\vartheta^\ast)}{\partial \vartheta}\sqrt{n}(\hat{\vartheta}_n - \vartheta_0) + \mathrm{o}_{\mathbb P}(1),
\end{aligned}
\end{equation*}
for some $\vartheta_i^\ast$, $i=1,\dots,2q+p+2$ between $\hat{\vartheta}_n$ and $\vartheta_0$.
In view of \eqref{res4}, there exists a neighborhood $\mathcal{V}(\vartheta_0)$ of $\vartheta_0$ such that
\begin{equation*}
\mathbb{E} \sup\limits_{\vartheta\in\mathcal{V}(\vartheta_0)}\left\| \dfrac{\partial^2 s_{t-h}(\vartheta)s_t(\vartheta)}{\partial \vartheta\partial \vartheta'}\right\| < \infty.
\end{equation*}
For a fixed $r_h$, using these inequalities, \eqref{res3} and Assumption \textbf{A0} ($\kappa_\eta <\infty$),
the almost sure convergence of $\vartheta^\ast$ to $\vartheta_0$, a second Taylor expansion and the ergodic theorem, we obtain
\begin{equation*}
\begin{aligned}
\dfrac{\partial r_h(\vartheta^\ast)}{\partial \vartheta} =\dfrac{\partial r_h(\vartheta_0)}{\partial \vartheta} + \mathrm{o}_{\mathbb P}(1) \underset{n\to \infty}{\longrightarrow} c_h := \mathbb{E}\left[s_{t-h}(\vartheta_0) \dfrac{\partial s_t(\vartheta_0)}{\partial \vartheta}\right] = -\mathbb{E}\left[s_{t-h}\dfrac{1}{\zeta_t^2(\vartheta_0)} \dfrac{\partial \zeta_t^2(\vartheta_0)}{\partial \vartheta}\right]
\end{aligned}
\end{equation*}
by the fact $\mathbb{E}[s_t(\vartheta_0)\partial s_{t-h}(\vartheta_0)/\partial \vartheta] = 0$. Note that, $c_h$ is the almost sure limit of the row $h$ of the matrix
$\hat{C}_m$. Consequently we have
\begin{equation}\label{C_m}
\dfrac{\partial \mathbf{r}_m(\vartheta_0)}{\partial \vartheta} \underset{n\to \infty}{\longrightarrow} C_m := \begin{pmatrix} c_1' \\ \vdots \\ c_m' \end{pmatrix}.
\end{equation}
It follows that
\begin{equation}\label{12}
\sqrt{n} \mathbf{\hat{r}}_m = \sqrt{n}\mathbf{r}_m + C_m\sqrt{n}(\hat{\vartheta}_n - \vartheta_0) + \mathrm{o}_{\mathbb P}(1).
\end{equation}
Denote $\sqrt{n}\mathbf{r}_m = n^{-1/2}\sum_{t=1}^n s_t\mathbf{s}_{t-1:t-m}$, where $\mathbf{s}_{t-1:t-m} = (s_{t-1},\ldots, s_{t-m})'$. We now derive the asymptotic distribution of $\sqrt{n}(\hat{\vartheta}'_n - \vartheta'_0, \mathbf{r}'_m)'$. In view of \eqref{res5}, the central limit theorem of \cite{Billingsley} applied to the martingale difference process
\begin{equation*}
\left\{\Upsilon_t = \left(s_t\dfrac{1}{\zeta_t^2(\vartheta_0)}\dfrac{\partial \zeta_t^2(\vartheta_0)}{\partial \vartheta'}, s_t\mathbf{s}_{t-1:t-m}'\right)' ; \sigma(\eta_u, u \leq t)\right\},
\end{equation*}
shows that
\begin{equation}\label{13}
\sqrt{n}(\hat{\vartheta}'_n - \vartheta'_0, \mathbf{r}'_m)'=\dfrac{1}{\sqrt{n}}\sum\limits_{t=1}^n\Upsilon_t  + \mathrm{o}_{\mathbb P}(1)\underset{n\to+\infty}{\overset{\mathcal{L}}{\longrightarrow}}  \mathcal{N}\left(0, \mathbb{E}[\Upsilon_t\Upsilon_t']\right),
\end{equation}
where
\begin{equation*}
\mathbb{E}\left[\Upsilon_t \Upsilon_t'\right] = (\kappa_\eta -1 ) \begin{pmatrix} J^{-1} & -J^{-1}C_m' \\ -C_mJ^{-1} & (\kappa_\eta - 1)I_m \end{pmatrix}.
\end{equation*}
Using \eqref{12} and \eqref{13} we obtain the distribution of $\sqrt{n}\mathbf{\hat{r}}_m$. Indeed $\sqrt{n}\mathbf{\hat{r}}_m \overset{\mathcal{L}}{\longrightarrow} \mathcal{N}(0, D)$ where $D$ is defined by
\begin{equation*}
D := (\kappa_\eta - 1)^2I_m - (\kappa_\eta - 1)C_m J^{-1}C_m'.
\end{equation*}

\textbf{$(iii)$ Invertibility of the matrix $D$.}

We now show that $D$ is invertible. Assumption \textbf{A5} entails that the law of $\eta_t^2$ is non degenerated, therefore $\kappa_\eta > 1$.
Thus study the invertibility of the matrix $D$ is similar to study the invertibility of $(\kappa_\eta - 1)I_m - C_m J^{-1} C_m'$.
Let
\begin{equation*}
V = \mathbf{s}_{t-1:t-m} + C_m J^{-1}\dfrac{1}{\zeta_t^2(\vartheta_0)} \dfrac{\partial \zeta_t^2(\vartheta_0)}{\partial \vartheta}\quad  \text{such that}\quad  \mathbb{E}\left[VV'\right] = (\kappa_\eta - 1)I_m - C_mJ^{-1}C_m'.
\end{equation*}
If the matrix $\mathbb{E}\left[VV'\right]$ is singular, then there exist a vector $\lambda = (\lambda_1,\ldots, \lambda_m)'$ not equal to zero such that
\begin{equation}\label{14}
\lambda' V = \lambda' \mathbf{s}_{t-1:t-m} + \lambda'C_mJ^{-1} \left(\dfrac{1}{\zeta_t^2(\vartheta_0)}\dfrac{\partial \zeta_t^2(\vartheta_0)}{\partial \theta}
+\dfrac{1}{\zeta_t^2(\vartheta_0)}\dfrac{\partial \zeta_t^2(\vartheta_0)}{\partial \tau}\right) = 0, \quad a.s.
\end{equation}
since $\vartheta=(\theta',\tau)'$.
Using the fact that
\begin{equation*}
\dfrac{1}{\zeta_t^2(\vartheta_0)}\dfrac{\partial \zeta_t^2(\vartheta_0)}{\partial \theta} = \dfrac{2}{\tau} \dfrac{1}{\zeta_t^\tau(\vartheta_0)}\dfrac{\partial \zeta_t^\tau(\vartheta_0)}{\partial \theta}\text{ and }
\dfrac{1}{\zeta_t^2(\vartheta_0)}\dfrac{\partial \zeta_t^2(\vartheta_0)}{\partial \tau} = -\dfrac{2}{\tau^2}\log(\zeta_t^\tau(\vartheta_0))+\dfrac{2}{\tau} \dfrac{1}{\zeta_t^\tau(\vartheta_0)}\dfrac{\partial \zeta_t^\tau(\vartheta_0)}{\partial \tau},
\end{equation*}
we can rewrite the equation \eqref{14} as follow
\begin{equation}\label{15}
\lambda' V = \lambda' \mathbf{s}_{t-1:t-m} + \mu' \dfrac{1}{\zeta_t^\tau(\vartheta_0)}\left(\tau \dfrac{\partial \zeta_t^\tau(\vartheta_0)}{\partial \theta}-\zeta_t^\tau(\vartheta_0)\log(\zeta_t^\tau(\vartheta_0))+\tau \dfrac{\partial \zeta_t^\tau(\vartheta_0)}{\partial \tau}\right)= 0, \quad a.s.
\end{equation}
with $\mu' = (2/\tau^2)\lambda'C_mJ^{-1}$. We remark that $\mu\neq 0$. Otherwise $\lambda' \mathbf{s}_{t-1:t-m}=0$
a.s., which implies that there exists $j\in\{1,...,m\}$ such that $s_{t-j}$ is measurable with
respect to the $\sigma-$field generated by $s_r$ for $t-1\le r\le t-m$ with $r\neq t-j$. This is impossible because the $s_t$'s are independent and non
degenerated.

%We can rewrite \eqref{15} such that the derivatives of $\zeta_t$ with respect to the vector $\theta$ and the unknown power $\tau$ appear in the expression.
We denote $\mu = (\nu_1', \nu_2)'$, where $\nu_1' = (\mu_1,\ldots, \mu_{2q+p+1})'$ and $\nu_2 = \mu_{2q+p+2}$; and we rewrite \eqref{15} as
\begin{equation*}
\lambda' V = \lambda' \mathbf{s}_{t-1:t-m} + \nu_1'\tau \dfrac{1}{\zeta_t^\tau(\vartheta_0)}\dfrac{\partial \zeta_t^\tau(\vartheta_0)}{\partial \theta} + \nu_2 \dfrac{1}{\zeta_t^\tau(\vartheta_0)}\left(-\zeta_t^\tau(\vartheta_0)\log(\zeta_t^\tau(\vartheta_0))+\tau \dfrac{\partial \zeta_t^\tau(\vartheta_0)}{\partial \tau}\right) = 0, \quad a.s.
\end{equation*}
or equivalent,
\begin{equation}\label{16}
\lambda' \mathbf{s}_{t-1:t-m}\zeta_t^\tau(\vartheta_0) + \nu_1'\tau \dfrac{\partial \zeta_t^\tau(\vartheta_0)}{\partial \theta} + \nu_2 \left(-\zeta_t^\tau(\vartheta_0)\log(\zeta_t^\tau(\vartheta_0))+\tau \dfrac{\partial \zeta_t^\tau(\vartheta_0)}{\partial \tau}\right) = 0, \quad a.s.
\end{equation}
The derivatives involved in \eqref{16} are defined recursively by
\begin{equation*}
\begin{aligned}
\dfrac{\partial \zeta_t^\tau(\vartheta)}{\partial \theta} &= \underline{{c}}_t(\vartheta) + \sum\limits_{j=1}^p \beta_j \dfrac{\partial \zeta_{t-j}^\tau(\vartheta)}{\partial \theta},\\
\dfrac{\partial \zeta_t^\tau(\vartheta)}{\partial \tau} &= \sum\limits_{i=1}^q \alpha_i^+\log(\varepsilon_{t-i}^+)(\varepsilon_{t-i}^+)^\tau + \alpha_i^-\log(-\varepsilon_{t-i}^-)(-\varepsilon_{t-i}^-)^\tau  + \sum\limits_{j=1}^p \beta_j \dfrac{\partial \zeta_{t-j}^\tau(\vartheta)}{\partial \tau},
\end{aligned}
\end{equation*}
where $c_t(\vartheta)$ is defined by replacing $\tilde\zeta^\tau_t(\vartheta)$ by $\zeta^\tau_t(\vartheta)$ in $\tilde c_t(\vartheta)$ (see \eqref{ctilde}).
We remind that  $\varepsilon_t^+ = \zeta_t\eta_t^+$ and $\varepsilon_t^- = \zeta_t\eta_t^-$ and let $R_t$ a random variable measurable with respect to $\sigma\{\eta_u, u \leq t\}$. We decompose \eqref{16} in four terms  and we have
\begin{align*}
\nu'_1\tau\dfrac{\partial \zeta_t^\tau(\vartheta_0)}{\partial \theta} & = \mu_2\tau\zeta_{t-1}^\tau(\eta_{t-1}^+)^\tau + \mu_{q+2}\tau\zeta_{t-1}^\tau(-\eta_{t-1}^-)^\tau + R_{t-2}, \\
\zeta_t^\tau &= \alpha_1^+\zeta_{t-1}^\tau(\eta_{t-1}^+)^\tau + \alpha_1^-\zeta_{t-1}^\tau(-\eta_{t-1}^-)^\tau + R_{t-2},\\
-\nu_2\zeta_t^\tau(\vartheta_0)\log(\zeta_t^\tau(\vartheta_0))&=-\nu_2\left(\alpha_1^+\zeta_{t-1}^\tau(\eta_{t-1}^+)^\tau + \alpha_1^-\zeta_{t-1}^\tau(-\eta_{t-1}^-)^\tau + R_{t-2}\right)
\\ &\quad \times\log\left(\alpha_1^+\zeta_{t-1}^\tau(\eta_{t-1}^+)^\tau + \alpha_1^-\zeta_{t-1}^\tau(-\eta_{t-1}^-)^\tau + R_{t-2}\right)
\\
%\nu_2\tau\dfrac{\partial \zeta_t^\tau(\vartheta_0)}{\partial \tau} & = \nu_2\tau\alpha_1^+\log\left(\zeta_{t-1}(\eta_{t-1}^+)\right)\zeta_{t-1}^\tau(\eta_{t-1}^+)^\tau + \nu_2\tau\alpha_1^-\log\left(\zeta_{t-1}(-\eta_{t-1}^-)\right)\zeta_{t-1}^\tau(-\eta_{t-1}^-)^\tau + R_{t-2},
%\\& = \nu_2\alpha_1^+\log\left(\zeta_{t-1}^\tau(\eta_{t-1}^+)^\tau\right)\zeta_{t-1}^\tau(\eta_{t-1}^+)^\tau + \nu_2
%\alpha_1^-\log\left(\zeta_{t-1}^\tau(-\eta_{t-1}^-)^\tau\right)\zeta_{t-1}^\tau(-\eta_{t-1}^-)^\tau + R_{t-2}, \\
\lambda' s_{t-1:t-m} &= \lambda_1 \eta_{t-1}^2 + R_{t-2},
\end{align*}
that gives
\begin{equation*}
\lambda' s_{t-1:t-m} \zeta_t^\tau = \lambda_1\zeta_{t-1}^\tau \left[ \alpha_1^+(\eta_{t-1}^+)^{\tau+2} + \alpha_1^-(-\eta_{t-1}^-)^{\tau + 2}\right] + \lambda_1\eta_{t-1}^2R_{t-2} + R_{t-2} + \left[(\eta_{t-1}^+)^\tau + (-\eta_{t-1}^-)^\tau\right]R_{t-2},
\end{equation*}
and
\begin{eqnarray*}
\nu_2\tau\dfrac{\partial \zeta_t^\tau(\vartheta_0)}{\partial \tau} & =& \nu_2\tau\alpha_1^+\log\left(\zeta_{t-1}(\eta_{t-1}^+)\right)\zeta_{t-1}^\tau(\eta_{t-1}^+)^\tau + \nu_2\tau\alpha_1^-\log\left(\zeta_{t-1}(-\eta_{t-1}^-)\right)\zeta_{t-1}^\tau(-\eta_{t-1}^-)^\tau + R_{t-2},
\\& =& \nu_2\alpha_1^+\log\left(\zeta_{t-1}^\tau(\eta_{t-1}^+)^\tau\right)\zeta_{t-1}^\tau(\eta_{t-1}^+)^\tau + \nu_2
\alpha_1^-\log\left(\zeta_{t-1}^\tau(-\eta_{t-1}^-)^\tau\right)\zeta_{t-1}^\tau(-\eta_{t-1}^-)^\tau + R_{t-2}.
\end{eqnarray*}
Following these previous expressions, \eqref{15} entails that almost surely
\begin{equation*}
\begin{aligned}%\label{final}
\lambda'V &= \lambda_1\zeta_{t-1}^\tau \left[ \alpha_1^+(\eta_{t-1}^+)^{\tau+2} + \alpha_1^-(-\eta_{t-1}^-)^{\tau + 2}\right] + \eta_{t-1}^2R_{t-2} + \left[R_{t-2}+\nu_2\alpha_1^+R_{t-2}\log(\zeta_{t-1}(\eta_{t-1}^+))\right](\eta_{t-1}^+)^\tau  \\
&\quad  + \left[R_{t-2}+\nu_2\alpha_1^-R_{t-2}\log(\zeta_{t-1}(-\eta_{t-1}^-))\right](-\eta_{t-1}^-)^\tau R_{t-2} + R_{t-2} \\
& -\nu_2\left(\alpha_1^+\zeta_{t-1}^\tau(\eta_{t-1}^+)^\tau + \alpha_1^-\zeta_{t-1}^\tau(-\eta_{t-1}^-)^\tau + R_{t-2}\right)
\log\left(\alpha_1^+\zeta_{t-1}^\tau(\eta_{t-1}^+)^\tau + \alpha_1^-\zeta_{t-1}^\tau(-\eta_{t-1}^-)^\tau + R_{t-2}\right)=0,
\end{aligned}
\end{equation*}
or equivalent to the two equations
\begin{equation}
\begin{aligned}
\label{17}
  \lambda_1\zeta_{t-1}^\tau \alpha_1^+(\eta_{t-1}^+)^{\tau+2} &-\left(\nu_2\alpha_1^+\zeta_{t-1}^\tau(\eta_{t-1}^+)^\tau + R_{t-2}\right)
\log\left(\alpha_1^+\zeta_{t-1}^\tau(\eta_{t-1}^+)^\tau + R_{t-2}\right)
 \\& +\left[R_{t-2} + \nu_2\alpha_1^+R_{t-2}\log(\zeta_{t-1}(\eta_{t-1}^+))\right](\eta_{t-1}^+)^\tau  + \eta_{t-1}^2R_{t-2}  + R_{t-2} =0,\, a.s.
\end{aligned}
\end{equation}
\begin{equation}
\begin{aligned}
\label{18}
\lambda_1\zeta_{t-1}^\tau \alpha_1^-(-\eta_{t-1}^+)^{\tau+2} &-\left(\nu_2\alpha_1^-\zeta_{t-1}^\tau(-\eta_{t-1}^-)^\tau + R_{t-2}\right)
\log\left(\alpha_1^-\zeta_{t-1}^\tau(-\eta_{t-1}^-)^\tau + R_{t-2}\right)
\\&+
\left[R_{t-2} + \nu_2\alpha_1^-R_{t-2}\log(\zeta_{t-1}(-\eta_{t-1}^-))\right](-\eta_{t-1}^-)^\tau   + \eta_{t-1}^2R_{t-2}  + R_{t-2} =0,\, a.s..
\end{aligned}
\end{equation}
%\begin{align}\label{17}
%  \lambda_1\zeta_{t-1}^\tau \alpha_1^+(\eta_{t-1}^+)^{\tau+2} &-\nu_2\left(\alpha_1^+\zeta_{t-1}^\tau(\eta_{t-1}^+)^\tau + R_{t-2}\right)
%\log\left(\alpha_1^+\zeta_{t-1}^\tau(\eta_{t-1}^+)^\tau + R_{t-2}\right)
%  +\left[R_{t-2} + \nu_2\alpha_1^+R_{t-2}\log(\zeta_{t-1}(\eta_{t-1}^+))\right](\eta_{t-1}^+)^\tau  + \eta_{t-1}^2R_{t-2}  + R_{t-2} =0 \\
%\label{18}
%\lambda_1\zeta_{t-1}^\tau \alpha_1^-(-\eta_{t-1}^+)^{\tau+2} &-\nu_2\left(\alpha_1^-\zeta_{t-1}^\tau(-\eta_{t-1}^-)^\tau + R_{t-2}\right)
%\log\left(\alpha_1^-\zeta_{t-1}^\tau(-\eta_{t-1}^-)^\tau + R_{t-2}\right)+
%\left[R_{t-2} + \nu_2\alpha_1^-R_{t-2}\log(\zeta_{t-1}(-\eta_{t-1}^-))\right](-\eta_{t-1}^-)^\tau   + \eta_{t-1}^2R_{t-2}  + R_{t-2} =0.
%\end{align}
Note that an equation of the form
\begin{equation*}
a\vert x \vert^{\tau + 2} + [b + c(\vert x \vert^\tau)]\log[b + c(\vert x \vert^\tau)]+[d + e\log(\vert x \vert)] \vert x \vert^\tau + f x^2 + g = 0
\end{equation*}
cannot have more than 11 positive roots or more than 11 negative roots, except if $a=b=c=d=e=f=g=0$. By assumption \textbf{A1}, Equations \eqref{17} and \eqref{18} thus imply that $\lambda_1(\alpha_1^+ + \alpha_1^-)=0$ and $\nu_2(\alpha_1^+ + \alpha_1^-)=0$.
If $\lambda_1=0$ and $\nu_2=0$ then $\lambda'\mathbf{s}_{t-1:t-m} := \lambda_{2:m}'\mathbf{s}_{t-2:t-m}$. By \eqref{16}, we can write that
\begin{equation*}
\left[\alpha_1^+\zeta_{t-1}^\tau(\eta_{t-1}^+)^\tau + \alpha_1^-\zeta_{t-1}^\tau(-\eta_{t-1}^-)^\tau \right]\lambda_{2:m}' \mathbf{s}_{t-2:t-m} = -\mu_2\zeta_{t-1}^\tau(\eta_{t-1}^+)^\tau + \mu_{q+2}\zeta_{t-1}^\tau(-\eta_{t-1}^-)^\tau + R_{t-2},
\end{equation*}
which entails
\begin{equation*}
\alpha_1^+\zeta_{t-1}^\tau(\eta_{t-1}^+)^\tau \lambda_{2:m}' \mathbf{s}_{t-2:t-m} = -\mu_2\zeta_{t-1}^\tau(\eta_{t-1}^+)^\tau + R_{t-2}
\end{equation*}
and a similar expression with $(-\eta_{t-1}^-)^\tau$ can be obtained. Subtracting the conditional expectation with respect to $\mathcal F_{t-2}=\sigma\{\eta_r^+, \eta_r^- \ {;}\  r\leq t-2\}$
in both sides of the previous equation, we obtain
\begin{equation*}
\begin{aligned}
\alpha_1^+ \zeta_{t-1}^\tau \lambda_{2:m}'\mathbf{s}_{t-2:t-m} \left[(\eta_{t-1}^+)^\tau - \mathbb{E}[(\eta_{t-1}^+)^\tau \vert \mathcal{F}_{t-2}]\right] &= \mu_2\zeta_{t-1}^\tau \left[\mathbb{E}[(\eta_{t-1}^+)^\tau \vert \mathcal{F}_{t-2}] - (\eta_{t-1}^+)^\tau\right],\quad a.s.\\
\alpha_1^+ \zeta_{t-1}^\tau \lambda_{2:m}'\mathbf{s}_{t-2:t-m} \left[(\eta_{t-1}^+)^\tau - \mathbb{E}[(\eta_{t-1}^+)^\tau]\right] &= \mu_2\zeta_{t-1}^\tau \left[\mathbb{E}[(\eta_{t-1}^+)^\tau] - (\eta_{t-1}^+)^\tau\right],\quad a.s..
\end{aligned}
\end{equation*}
Since the law of $\eta_t$ is non degenerated, we have $\alpha_1^+ = \mu_2 = 0$ and symmetrically $\alpha_1^- = \mu_{q+2} = 0$. But for APGARCH$(p,1)$ models, it is impossible to have $\alpha_1^+ = \alpha_1^- = 0$ by the assumption \textbf{A4}. The invertibility of $D$ is thus shown in this case.
For APGARCH$(p,q)$ models, by iterating the previous arguments, we can show by induction that \eqref{15} entails $\alpha_1^+ + \alpha_1^- = \ldots = \alpha_q^+ + \alpha_q^-=0$.
Thus $\lambda_1=\dots=\lambda_m=0$ which leads to a contradiction. The non-singularity of the matrix $D$ follows.\zak
\subsection{Proof of Theorem~\ref{resultat}}
The almost sure convergence of $\hat{D}$ to $D$ as $n$ goes to infinity is easy to show using the consistency result.
%and is omitted.
The matrix $D$ can be rewritten as
$D = (\kappa_\eta - \hat{\kappa}_\eta)B + (\hat{\kappa}_\eta-1)A,$
where %$A = C_mJ^{-1}C_m'$ and $B = (\kappa_\eta - 1)I_m - C_mJ^{-1}C_m$.
the matrices $A$ and $B$ are given by
\begin{align*}
A%&=(C_m - \hat{C}_m)J^{-1}C_m' + \hat{C}_m(J^{-1} - \hat{J}^{-1})C_m' + \hat{C}_m\hat{J}^{-1}(C_m' - \hat{C}_m') + \hat{C}_m\hat{J}^{-1}\hat{C}_m',\\
&= (C_m - \hat{C}_m)J^{-1}C_m' + \hat{C}_m(J^{-1} - \hat{J}^{-1})C_m' + \hat{C}_m\hat{J}^{-1}(C_m' - \hat{C}_m') + \hat{A},\\
B%&=(\kappa_\eta - \hat{\kappa}_\eta)I_m + (\hat{\kappa}_\eta - 1)I_m + (A-\hat{A}) + \hat{A},\\
&= (A-\hat{A}) + (\kappa_\eta - \hat{\kappa}_\eta)I_m + \hat{B},
\end{align*}
with
$\hat{A} = \hat{C}_m\hat{J}^{-1}\hat{C}_m'$ and $\hat{B}=(\hat{\kappa}_\eta-1)I_m - \hat{A}$.
Finally, we have
$$D-\hat{D} = (\kappa_\eta - \hat{\kappa}_\eta)B + (\hat{\kappa}_\eta - 1)\left[(A-\hat{A}) + (\kappa_\eta - \hat{\kappa}_\eta)I_m\right].$$
For any multiplicative norm, we have
$$\Vert D-\hat{D}\Vert \leq \vert \kappa_\eta - \hat{\kappa}_\eta \vert \Vert B \Vert + \vert \hat{\kappa}_\eta - 1 \vert \left[\Vert A - \hat{A}\Vert + \vert\kappa_\eta - \hat{\kappa}_\eta \vert m \right]$$
and
%\begin{align*}
$$
\Vert A  - \hat{A} \Vert %&\leq \Vert C_m - \hat{C}_m \Vert \Vert J^{-1} \Vert \Vert C_m'\Vert + \Vert \hat{C}_m\Vert \Vert J^{-1} - \hat{J}^{-1} \Vert \Vert C_m' \Vert + \Vert C_m \Vert \Vert \hat{J}^{-1} \Vert \Vert C_m' - \hat{C}_m'\Vert \\
\leq \Vert C_m - \hat{C}_m \Vert \Vert J^{-1} \Vert \Vert C_m' \Vert  + \Vert \hat{C}_m\Vert \Vert J^{-1} \Vert \Vert \hat{J} - J \Vert \Vert \hat{J}^{-1} \Vert \Vert C_m' \Vert + \Vert C_m \Vert \Vert \hat{J}^{-1}\Vert \Vert C_m' - \hat{C}_m'\Vert. $$
%\\&\underset{n\to+\infty}{\longrightarrow} 0, a.s.
%\end{align*}
In view of \eqref{res3}, we have $\Vert C_m \Vert<\infty$. Because the matrix $J$ is nonsingular, we have $\Vert J^{-1} \Vert<\infty$ and $$\Vert \hat{J}^{-1}-J^{-1} \Vert\underset{n\to+\infty}{\longrightarrow} 0,\quad a.s.$$
by consistency of $\hat\vartheta_n$. Under
Assumption \textbf{A5}, we have $\vert \kappa_\eta - 1\vert \leq K$. Using the previous arguments and also the strong consistency of  $\hat\vartheta_n$, we have
$$\vert \kappa_\eta - \hat{\kappa}_\eta\vert\underset{n\to+\infty}{\longrightarrow} 0,\quad a.s.\text{ and } \Vert C_m - \hat{C}_m \Vert\underset{n\to+\infty}{\longrightarrow} 0,\quad a.s.$$ We then deduce that $\Vert B \Vert \leq K$ and the conclusion follows. Thus $\hat{D} \to D$ almost surely, when $n \to +\infty$.

To conclude the proof of Theorem~\ref{resultat}, it suffices to use Theorem~\ref{resultat.r_h} and the following result:
if $\sqrt{n} \mathbf{\hat{r}}_m \overset{\mathcal{L}}{\longrightarrow} \mathcal{N}\left(0,D\right)$, with $D$
nonsingular, and if $\hat{D}\to D$ in probability, then $
n \mathbf{\hat{r}}_m'\hat{D}^{-1}\mathbf{\hat{r}}_m \overset{\mathcal{L}}{\longrightarrow} \chi^2_m.$
\zak
\subsection{Condition of strict stationarity of model \eqref{model}}
The probabilistic properties of the model \eqref{model} rely on the sequence of matrices $(C_{0t})$ defined by
\begin{equation}\label{cot}
C_{0t} = \begin{pmatrix} \kappa(\eta_t) & \beta_{0p} & \alpha_{[2:q-1]} & \alpha_{[q:q]} \\
I_{p-1} & 0_{(p-1)\times 1} & 0_{(p-1)\times 2(q-2)}  & 0_{(p-1)\times 2} \\
\underline{\eta}_t & 0_{2\times 1} & 0_{2\times 2(q-2)} & 0_{2\times 2} \\
0_{2(q-2)\times (p-1)} & 0_{2(q-2)\times 1} & I_{2(q-2)} & 0_{2(q-2)\times 2} \end{pmatrix},
\end{equation}
where $I_k$ denotes the identity matrix of size $k$ and, for $i \leq j$,
\begin{equation*}
\begin{aligned}
&\kappa(\eta_t) = \big (\beta_{01}+\alpha_{01}^+(\eta_t^+)^\delta + \alpha_{01}^-(-\eta_t^-)^\delta ,\beta_{02},\ldots, \beta_{0p-1}\big ),\\
&\alpha_{[i:j]} = (\alpha_{0i}^+,\alpha_{0i}^-,\ldots, \alpha_{0j}^+, \alpha_{0j}^-),\qquad \underline{\eta}_t = \begin{pmatrix} (\eta_t^+)^\delta & 0_{1\times (p-1)} \\
(-\eta_{t}^-)^\delta & 0_{1\times(p-2)}\end{pmatrix}.
\end{aligned}
\end{equation*}
%%%% Debut Bibliographie

\bibliographystyle{apalike}
\bibliography{biblio-oth}

%%%% Fin Bibliographie

\end{document}